\documentclass[11pt]{article}
\usepackage{amssymb,amsfonts,amsmath,latexsym,epsf,tikz,url}

\newtheorem{theorem}{Theorem}[section]

\newtheorem{corollary}[theorem]{Corollary}
\newtheorem{lemma}[theorem]{Lemma}

\newtheorem{example}[theorem]{Example}
\newcommand{\proof}{\noindent{\bf Proof.\ }}
\newcommand{\qed}{\hfill $\square$\medskip}

\begin{document}

\title{Distinguishing number and distinguishing index of some operations  on graphs}

\author{
Saeid Alikhani  $^{}$\footnote{Corresponding author}
\and
Samaneh Soltani
}

\date{\today}

\maketitle

\begin{center}
Department of Mathematics, Yazd University, 89195-741, Yazd, Iran\\
{\tt alikhani@yazd.ac.ir, s.soltani1979@gmail.com}
\end{center}

\begin{abstract}
The distinguishing number (index) $D(G)$ ($D'(G)$) of a graph $G$ is the least integer $d$
such that $G$ has an vertex labeling (edge labeling)  with $d$ labels  that is preserved only by a trivial automorphism.  We examine the effects on $D(G)$ and $D'(G)$ when $G$ is modified by operations on vertex and  edge of $G$. Let $G$ be a connected graph of order $n\geq 3$. We show that $-1\leq D(G-v)-D(G)\leq D(G)$, where $G-v$ denotes the graph obtained from $G$ by removal of a vertex  $v$ and all edges incident to $v$ and these inequalities are true for the distinguishing index. Also we prove that $|D(G-e)-D(G)|\leq 2$ and $-1 \leq D'(G-e)-D'(G)\leq 2$, where $G-e$ denotes the graph obtained from $G$ by simply removing the edge $e$.  Finally we consider the vertex contraction and the edge contraction of $G$ and prove that the edge contraction decrease the distinguishing number (index) of $G$ by at most one and increase by at most  $3D(G)$ ($3D'(G)$). 
\end{abstract}

\noindent{\bf Keywords:}  Distinguishing index; distinguishing number; edge contraction.  

\medskip
\noindent{\bf AMS Subj.\ Class.:} 05C15, 05E18

\section{Introduction}

Let $G = (V ,E)$ be a simple  graph with $n$ vertices.  We use the standard graph notation (\cite{Sandi}).  
The set of all automorphisms of $G$, with the operation of composition of permutations, is a permutation group
on $V$ and is denoted by $Aut(G)$.  
A labeling of $G$, $\phi : V \rightarrow \{1, 2, \ldots , r\}$, is  $r$-distinguishing, 
if no non-trivial  automorphism of $G$ preserves all of the vertex labels.
In other words,  $\phi$ is $r$-distinguishing if for every non-trivial $\sigma \in Aut(G)$, there
exists $x$ in $V$ such that $\phi(x) \neq \phi(x\sigma)$. 
The distinguishing number of a graph $G$ has defined by Albertson and Collins \cite{Albert} and  is the minimum number $r$ such that $G$ has a labeling that is $r$-distinguishing.  
Similar to this definition, Kalinkowski and Pil\'sniak \cite{Kali1} have defined the distinguishing index $D'(G)$ of $G$ which is  the least integer $d$
such that $G$ has an edge colouring   with $d$ colours that is preserved only by a trivial
automorphism.  These indices  has developed  and number of papers published on this subject (see, for example \cite{soltani2,Klavzar,fish}).

We use the following notations: The set of vertices adjacent in $G$ to a vertex of a
vertex subset $W\subseteq  V$ is the open neighborhood $N_G(W )$ of $W$. The closed neighborhood  $G[W ]$ also includes all vertices of $W$ itself. In case of a singleton set $W =\{v\}$ we write $N_G(v)$ and $N_G[v]$ instead of $N_G(\{v\})$ and $N_G[\{v\}]$, respectively.  We omit the subscript when the graph $G$ is clear from the context.  The complement of  $N[v]$ in $V(G)$ is  denoted by $\overline{N[v]}$. 
The graph $G-v$ is a graph that is made  by deleting the vertex $v$ and all edges connected to $v$ from the graph $G$. Similarly, the graph $G-e$ is a graph that obtained from $G$ by simply removing the  edge $e$.  An operation on vertex $v$ of $G$ which is denoted by $G\odot v$  is a graph obtained by the removal of all edges between any pair of neighbors of $v$. The contraction of $v$ in $G$ denoted by $G\circ v$ is the graph obtained by deleting $v$ and putting a clique on the (open) neighbourhood of $v$. The  contraction of an edge $e$ with endpoints $u,v$ in graph $G$ is denoted by $G\circ e$ and is the replacement of $u$ and $v$ with a single vertex such that edges incident to the new vertex are the edges other than $e$ that were incident with $u$ or $v$.

In the study of distinguishing number and distinguishing index of graphs,  this naturally raises the question:  What happens to the distinguishing number and the  distinguishing index, when we consider some operations on the vertices and  the edges of a graph? 
In this paper we would like to answer to this question. 

\medskip 

In the next section, we examine the effects on $D(G)$ and $D'(G)$ when $G$ is modified by deleting a vertex or deleting an edge. In the Section 3, we consider the distinguishing number and the  distinguishing index of  $G\odot v$.  In the last section we  study the effects on $D(G)$ and $D'(G)$ when $G$ is modified by contracting a vertex and contracting an  edge.


\section{Vertex and edge removal }
The graph $G-v$ is a graph that is made  by deleting the vertex $v$ and all edges connected to $v$ from the graph $G$ and the graph $G-e$ is a graph that obtained from $G$ by simply removing the  edge $e$. Our main results in this section are  in obtaining  a  bound for distinguishing number and index of $G-v$ and $G-e$. To do this we need to consider some preliminaries.  A graph $G$ is almost  spanned by a subgraph $H$ if $G-v$ is spanned by $H$ for some $v \in V (G)$.

\begin{lemma}{\rm\cite{nord}}\label{lemm1}
If a graph $G$ is spanned or almost spanned by a subgraph $H$, then $D'(G) \leqslant D'(H) + 1$.
\end{lemma}
\proof
 We colour the edges of $H$ with colours $1,\ldots ,D'(H)$, and all other edges of $G$ with an additional colour $0$. If $\varphi$ is an automorphism of $G$ preserving
this colouring, then $\varphi(x) = x$, for each $x \in V (H)$. Moreover, if $H$
is a spanning subgraph of $G-v$, then also $\varphi(v) = v$. Therefore, $\varphi$ is the
identity  map 
and we have the result.\qed

The following result shows that removing a vertex of $G$ can decrease the distinguishing number (index) by at most one but can increase by at most to double  of distinguishing number (index) of $G$.  
\begin{theorem}\label{thmG-v}
Let $G$ be a connected graph of order $n\geqslant 3$ and $v\in V(G)$. Then we have
\begin{itemize}
\item[(i)] $D(G) -1\leqslant D(G-v)\leqslant 2D(G)$.
\item[(ii)] $D'(G) -1\leqslant D'(G-v)\leqslant 2D'(G)$.
\end{itemize}
\end{theorem}
\proof 
\begin{enumerate}
\item[(i)]
 We label the vertices of $G-v$ with labels $1,\ldots, D(G-v)$, and the vertex $v$ of $G$ with an additional label $0$. If $f$ is an automorphism of $G$ preserving this labeling, then $f(v)=v$ and the restriction of $f$ to  $G-v$ is an automorphism of $G-v$ preserving this labeling. Since we labeled $G-v$ in a distinguishing way at first, so $f$ is the identity map. Thus we proved the first inequality in this part, similar to the proof of Lemma \ref{lemm1}. 

For the second inequality, let $N_G(v)=\{v_1,\ldots , v_m\}$. We label the vertices of $G$ with labels $1,\ldots , D(G)$, in a distinguishing way. Next we change the label of the elements of $N_G(v)$ such that if the label of $v_i$, $1\leqslant i \leqslant n$, is $t_i$, then we change this label to $t_i + D(G)$. This labeling is a distinguishing labeling for $G-v$. Because:

If  $f$ is an  automorphism of $G-v$ preserving the labeling, then $f(N_G(v))= N_G(v)$ and $f(\overline{N[v]})=\overline{N[v]}$, then by defining $f(v)=v$, $f$ is the automorphism of $G$ preserving the labeling, because $f$ preserves the adjacency relation on $G$. Regarding to the method of the labeling of $G-v$, the map  $f$ is the identity map on $G-v$, and so $f$ is the identity map on $G$.
Since we used at most $2D(G)$ labels, we have $D(G-v)\leqslant 2D(G)$.

\item[(ii)]  The first inequality follows directly  from Lemma \ref{lemm1}. For the second inequality, let $N_G(v)=\{v_1,\ldots , v_m\}$ . We label the edges of $G$ with labels $1,\ldots , D'(G)$, in a distinguishing way. Next we change the label of the edges $\{vv_1,\ldots , vv_m\}$ such that if the label of $vv_i$, $1\leqslant i \leqslant n$, is $t_i$, then we change this label to $t_i + D'(G)$. The rest of proof is similar to the proof of Part (i). \qed
\end{enumerate}

\begin{figure}[ht]
\begin{center}
\includegraphics[width=0.4\textwidth]{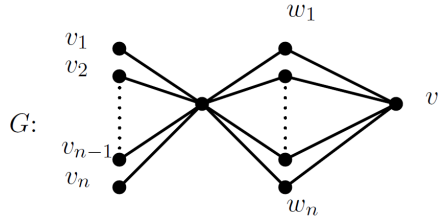}
	\caption{\label{fig6}The graph $G$ with $D(G-v)=2D(G)=2n$.}
\end{center}
\end{figure}

\begin{example}\label{examp2}
In this example we present some graphs such that they obtain the bounds of Theorem \ref{thmG-v}. Note that these examples can be used for distinguishing index, too. 

\begin{itemize}
\item[(i)] Let $K_{1,n}$ be the star graph such that $n\geqslant 3$ and let $v$ be an its arbitrary end vertex (vertex of degree one). Then $D(K_{1,n} -v)=D(K_{1,n-1})=n-1$, and so $D(K_{1,n})-1= D(K_{1,n} -v)$. 
\item[(ii)] Let $G$ be a graph has  shown  in Figure \ref{fig6}. With regard to the degree sequences  of $G$, it is easy to compute its automorphism group, and hence we can get $D(G)=n$. Also $D(G-v)=D(K_{1,2n})=2n$, and so $D(G-v)= 2D(G)$.
\end{itemize}
\end{example}

\begin{figure}[ht]
\begin{center}
\includegraphics[width=0.7\textwidth]{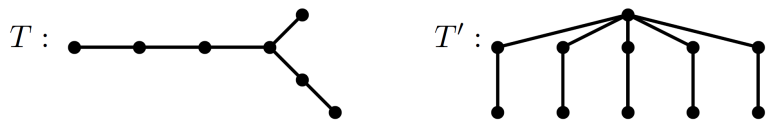}
	\caption{\label{fig7} Bicentered and a centered  tree, respectively.}
\end{center}
\end{figure}

For some parameters of a graph such as domination number $\gamma(G)$, has proved that for every tree $T$ there exists a vertex $v$ such that $\gamma(T-v)=\gamma(T)$ (\cite{domination}). Here we show that this is not true for the distinguishing number and the distinguishing index. 
Let $T$ be a tree of order $n\geqslant 3$.  We present  two tree (a centered tree  and bicentered  tree) such that $D(T-v)\neq D(T)$ and $D'(T-v)\neq D'(T)$ for all $v\in V(T)$. Let $T$ be as shown in Figure \ref{fig7}. Since $T$ is an asymmetric graph, so $D(T) =1$. It can be seen that $2= D(T-v)\neq 1$ for all $v\in V(T)$ by trial and error. For  a centered  graph, see the graph $T'$ in   Figure \ref{fig7}. This kind of tree some times called a spider. More precisely, a spider is the graph formed by subdividing all of the edges of a star $K_{1,t}$. 
 It is easy to obtain that  $D(T') =3$ and $D(T'-v)\neq 3$ for all $v\in V(T')$.  
 Note that the number of branches in $T'$ is five. Since the automorphisms of $T'$ are the permutations of branches,  we  have  the following theorem: 
 
 \begin{theorem} 
 Suppose that $T'$ is  a spider which has formed by subdividing all of the edges of a star $K_{1,t}$. If $t=k^2+1$, for some $k\geqslant 2$, then $D(T'-v)\neq D(T')$ for all $v\in V(T')$. 
 \end{theorem}

Now we  examine the effects on $D(G)$ and $D'(G)$ when $G$ is modified by deleting an edge.

\begin{theorem}\label{thmG-e}
For each connected graph  $G$ and $e\in E(G)$, $\vert D(G-e)-D(G) \vert \leqslant 2$.
\end{theorem}
\proof
We first prove the inequality    $D(G-e)\leqslant D(G)+2$.  We define a vertex distinguishing labeling with $D(G)+2$ labels for $G-e$. First we label the vertices of $G$ with $D(G)$ labels in a distinguishing way, next we replace the label of $v$ and $w$ by two new labels. This vertex labeling is distinguishing for $G-e$, because if $f$ is an automorphism of $G-e$ such that it preserves the labeling, then $f(v)= v$ and $f(w)=w$, and so $f$ is an automorphism of $G$. Since we labeled $G$ in a  distinguishing way at first, $f$ is the identity automorphism.

Now we shall prove that  $D(G)\leqslant D(G-e)+2$. For this case,  we label the vertices of  $G$ with $D(G-e)+2$ labels in a distinguishing way. First we label $G-e$ with $D(G-e)$ labels in a distinguishing way. Next we replace the label of $v$ and $w$ by two new labels. This vertex labeling is distinguishing for $G$, because if $f$ is an automorphism of $G$ such that it preserves the labeling, then $f(v)= v$ and $f(w)=w$, and so $f$ is an automorphism of $G-e$. Since we labeled $G-e$ in a  distinguishing way at first,  $f$ is the identity automorphism.
Therefore  we have the result.\qed

The bounds of Theorem \ref{thmG-e} are sharp. It is immediate  that for the complete graph of order $n$, $K_n$, $D(K_n)=n$, although a simple computation yields that $D(K_n-e)=n-2$. Let $G=K_1+(K_2\cup (n-2)K_1)$ and $e$ be an edge of $K_2$.  Since $G-e=K_{1,n}$, so  $D(G-e)=n$, even though  $D(G)=n-2$.

\medskip 
The following theorem examine the effect on $D'(G)$ when $G$ is modified by deleting an edge. 
\begin{theorem}
For each  connected graph  $G$ and $e\in E(G)$,  
$-1 \leqslant D'(G-e)-D'(G)\leqslant  2.$
\end{theorem}
\proof
We first prove the inequality $D'(G)-1 \leqslant D'(G-e)$.  We define an edge distinguishing labeling with $D'(G-e)+1$ labels for $G$. First we label the edges of $G-e$ with $D'(G-e)$ labels in a distinguishing way, next we assign a new label to the edge $e$.  This edge labeling is distinguishing for $G$. Because if $f$ is an automorphism of $G$ such that it preserves the labeling, then $f(\{v,w\})= \{v,w\}$, and so $f$ is an automorphism of $G-e$. Since we labeled $G-e$ in a  distinguishing way at first, the map  $f$ is the identity automorphism. This bound is sharp, because $D'(K_{1,n})=n$ and $D'(K_{1,n}-e)=n-1$ where $e$ is an arbitrary edge of $K_{1,n}$.

Now we shall prove that $D'(G-e)\leqslant D'(G) + 2$. For this step we label the edges of $G$ with $D'(G)$ labels in a distinguishing way.  Next we replace the label of all  incident edges  to $v$ except $e$, by a new label, and also we replace the label of all incident edges to $w$ except $e$, by another new label.
  This edge labeling is distinguishing for $G-e$, because if $f$ is an automorphism of $G-e$ such that it preserves the labeling, then $f(v)= v$ and $f(w)=w$, and so $f$ is an automorphism of $G$. Since we labeled the edges of $G$ in a  distinguishing way at first,  $f$ is the identity automorphism. This bound is sharp,  because it is sufficient to consider  $G=K_1+(K_2\cup (n-2)K_1)$ and $e$ be an edge of $K_2$.  Since $G-e=K_{1,n}$, so  $D'(G-e)=n$, even though  $D'(G)=n-2$. \qed

\section{Distinguishing number and distinguishing index of $G\odot v$}

 We denote by $G\odot v$ the graph obtained from G by the removal of all edges between
any pair of neighbors of $v$, note $v$ is not removed from the graph \cite{alikhani}. It is clear that if $G$ is a graph of order $n$ and the degree of a  vertex $v$ is $n-1$ in $G$, then $G\odot v$ is a star graph as $K_{1,n-1}$, and so $D(G\odot v) = n-1$ and $D'(G\odot v) = n-1$. Therefore  for every vertex $v$ of $G$,  $D(G\odot v) \leqslant n-1$ and $D'(G\odot v) \leqslant n-1$. The following theorem is a lower bound for distinguishing number and index of $G\odot v$.

\begin{theorem}\label{odot}
Let  $G$   be a connected graph of order $n$ and $v$ be a vertex of $G$. Then 
\begin{itemize}
\item[(i)] $D(G)-1\leqslant D(G\odot v)$,
\item[(ii)] $D'(G)-1\leqslant D'(G\odot v)$.
\end{itemize}

\end{theorem}
\proof
(i) Let $N(v)=\{x_1,\ldots , x_s,y_1,\ldots , y_t\}$ where $0\leqslant s,t \leqslant n-1$ and let the degree of $x_i$'s  in $G\odot v$ be  one and the degree of $y_j$'s in $G\odot v$ be  greater than one. So there exist the vertices  $z_j$'s,  $1\leqslant j \leqslant t$, such that the distance between $z_j$ and $y_j$ in $G$ is one and the distance between $z_j$ and $v$ in $G$ is two. First note that if $f$ is an automorphism of $G$  such that $f(v)=v$ then $f(\{x_1,\ldots , x_s\})=\{x_1,\ldots , x_s\}$ and $f(\{y_1,\ldots , y_t\})=\{y_1,\ldots , y_t\}$. Because if $f(v)=v$ and there exists $1\leqslant i\leqslant s$ and $1\leqslant j\leqslant t$ such that $f(y_j)=x_i$, then the distance between $f(z_j)$ and $x_i$ in $G$ should be one, and so $f(z_j)\in \{x_1,\ldots , x_s,y_1,\ldots , y_t\}$. On the other hand  the distance between $f(z_j)$ and $v$ in $G$ should be  two.  Since $f(z_j)\in \{x_1,\ldots , x_s,y_1,\ldots , y_t\}$ and the distance between each element  of this set and the vertex $v$ is one, we have a contradiction. Therefore if $f$ is an automorphism of $G$  such that $f(v)=v$ then the restriction of $f$ to $G\odot v$ is an automorphism of $G\odot v$.

Now we define a distinguishing labeling for $G$ with $D(G\odot v)+1$ labels. First we label $G\odot v$ with $D(G\odot v)$ labels in a distinguishing way. Next we replace the label of $v$ by a new label, and so the graph $G\odot v$ is labeled in a distinguishing way with $D(G\odot v)+1$ labels. This labeling is a distinguishing labeling for $G$. Because  if $f$ is an automorphism of $G$ such that it preserves the labeling, then $f(v)=v$. We have seen that the restriction of $f$ to $G\odot v$ is an automorphism of $G\odot v$. Since this labeling is a distinguishing labeling for $G\odot v$, so $f$ is the identity automorphism.

(ii) We define an edge distinguishing labeling for $G$ with $D'(G\odot v)+1$ labels. First we label the edges of  $G\odot v$ with $D'(G\odot v)$ labels in a distinguishing way. Next we assign all the removed edges of $G$, a new label. This edge labeling is distinguishing for $G$, because 
 if $f$ is an automorphism of $G$ such that it preserves the labeling, then $f$ maps all the removed edges to  the removed edges, and so the restriction of $f$ to $G\odot v$ is an automorphism of $G\odot v$. Since we labeled $G\odot v$ in a distinguishing way at first, so $f$ is  the identity automorphism. \qed

It is easy to see that the lower bounds in Theorem \ref{odot} are sharp. We know $D(K_n)=n$, although $D(K_n\odot v) =D(K_{1,n-1})=n-1$. Let the  friendship graph $F_n$ be  the join of $K_1$ with $n$-copies of $K_2$, i.e., $F_n=K_1+nK_2$. If  $v$ is a non-central vertex of $F_3$, then $D'(F_3\odot v)=2$, however $D'(F_3)=3$. By the following theorem which gives the distinguishing number and the distinguishing index of $F_n$, we can find numerous $n$ such that $D'(F_n\odot v)= D'(F_n)-1$ where $v$ is a non-central vertex of $F_n$.   
  
  \begin{theorem}\label{soltani}{\rm{\cite{soltani2}}}
  	The distinguishing number of the friendship graph $F_n$  $(n\geq 2)$ is  $$D(F_n)= \big\lfloor \frac{1+\sqrt{8n+1}}{2}\big\rfloor.$$ 
  \end{theorem}

We end this section with the following theorem: 
\begin{theorem}
	\begin{enumerate}
		\item [(i)]   There exists a graph $G$ and a vertex of $G$, such  that $\frac{D(G\odot v)}{D(G)}$ is arbitrarily large, 
		\item[(ii)]    There exists a graph $G$ and a vertex of $G$, such  that $\frac{D'(G\odot v)}{D'(G)}$ is arbitrarily large. 
		\end{enumerate}
\end{theorem}
\proof 
\begin{enumerate}
	\item [(i)] 
	Let $v$  be central vertex of the friendship graph  $F_n$.  Then  $$\frac{D(F_n\odot v)}{D(F_n)}=\frac{D(K_{1,2n})}{D(F_n)}=\frac{2n}{\lfloor \frac{1+\sqrt{8n+1}}{2}\rfloor},$$   
	and this fraction can be arbitrarily large. 
	
	\item[(ii)] Since $D'(K_n)=2$ for $n>3$, and $D'(K_n\odot v) =D'(K_{1,n-1})=n-1$, we can make $\frac{D'(G\odot v)}{D'(G)}$ arbitrarily large.\qed
\end{enumerate}

\section{Distinguishing number and distinguishing index of vertex and edge contraction}
 Let $v$ be a vertex in $G$. The contraction of $v$ in $G$ denoted by $G\circ v$ is the graph obtained by deleting $v$ and putting a clique on the (open) neighbourhood of $v$. Note that this operation does not create parallel edges; if two neighbours of $v$ are already adjacent, then they remain simply adjacent (see \cite{Walsh}). 

\begin{theorem}\label{vertexcontraction}
Let $G$ be a connected graph of order $n\geqslant 2$ and $v$ be a vertex of $G$.   Then
\begin{itemize}
\item[(i)] $D(G)-1\leqslant D(G\circ v)$,
\item[(ii)] $ D'(G\circ v)\leqslant D'(G)+1$.
\end{itemize}
\end{theorem}
\proof  (i) We define a distinguishing vertex labeling for $G$ with $D(G\circ v)+1$ labels. First we label $G\circ v$ with $D(G\circ v)$ labels in a distinguishing way. Assigning a new label to $v$,  we consider  this labeling, as a labeling of $G$. This labeling is a distinguishing vertex labeling for $G$, because if $f$ is an automorphism of $G$ preserving the labeling, then $f(v)=v$, and so $f(N_G(v))=N_G(v)$. Therefore $f$ is an automorphism of $G\circ v$. Since we labeled $G\circ v$ in a distinguishing way at first, $f$ is the identity automorphism.

(ii)  If the valency of the vertex $v$ of $G$ is full, then with respect to the degree of the vertex $v$, it can be seen that  $ D'(G\circ v)\leqslant D'(G)+1$. So we suppose that the valency of the vertex $v$ is not full and  label the edges of $G$ with the labels $\{1,\ldots , D'(G)\}$  in a distinguishing way. To continue the proof, we consider two following  cases:

Case 1) If $\vert N_G(v)\vert\neq 2$, then assign  the  labels $0$ and  $1$ to  the  edges of the clique $K_{N_G(v)}$  in construction $G\circ v$,  such that the edges of the clique $K_{N_G(v)}$, have been labeled in a distinguishing way. Now we consider  this labeling for $G\circ v$. Note that since the valency of the vertex $v$ is not full, we can label the edges of the clique $K_{N_G(v)}$ with two labels in a distinguishing way. This labeling is a distinguishing edge labeling for $G\circ v$, because if $f$ is an automorphism of $G\circ v$ preserving the labeling then $f(N_G(v))= N_G(v)$. Since the edges of the clique $K_{N_G(v)}$ have been labeled in a distinguishing way, so $f(x)=x$ for each $x\in N_G(v)$, and so by defining $f(v)=v$ we can extend $f$ to an automorphism of $G$. Since we labeled $G$ in a distinguishing way, $f$ is the identity automorphism.

 Case 2) If $\vert N_G(v)\vert= 2$, we can assume that $N_G(v)=\{w,w'\}$. Since the valency of the vertex $v$ is not full, so without loss of generality we suppose that there exists the new vertex  $z$ that is incident to $w$. Assigning the edges $zw$ and $ww'$ the new label $0$ and  $1$, respectively. We consider  this labeling for $G\circ v$.  This labeling is a distinguishing edge labeling for $G\circ v$, because if $f$ is an automorphism of $G\circ v$ preserving the labeling then $f(w)= w$ and $f(w')=w'$, and so by defining $f(v)=v$ we can extend $f$ to an automorphism of $G$. Since we labeled $G$ in a distinguishing way, $f$ is the identity automorphism.  \qed

\begin{figure}[ht]
\begin{center}
\includegraphics[width=0.9\textwidth]{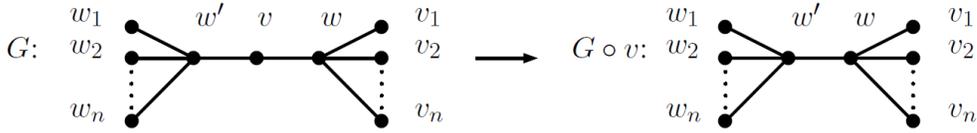}
	\caption{\label{fig9} $D'(G\circ v)=D'(G)+1$.}
\end{center}
\end{figure}

The bounds presented in  Theorem \ref{vertexcontraction} for the distinguishing number and index are sharp. For the Part (i), it is sufficient to consider $G=K_{1,n}$ and $v$  a vertex of degree one.   For the Part (ii),  let  $G$ be as shown in Figure \ref{fig9}. In regard to the degree of the vertices of $G$ we can obtain the automorphism group, and so we can get $D'(G)=n$ and $D'(G\circ v)=n+1$.

\begin{corollary}
\begin{itemize}
\item[(i)] There is a graph $G$ and a vertex $v$ of $G$ such that  the value of $\frac{D(G\circ  v)}{D(G)}$ can be  arbitrarily large. 
\item[(ii)]  There is a graph $G$ and a vertex $v$ of $G$ such that the value of $\frac{D'(G\circ  v)}{D'(G)}$ can be  arbitrarily small.
\end{itemize}
\end{corollary}
\proof (i)  We consider the friendship graph $F_n$. If $v$  is the  central vertex of $F_n$,  then $F_n\circ v$ is the complete graph $K_{2n}$, and so $D(F_n\circ  v)=2n$. Now  by Theorem \ref{soltani} it can be seen that the value $\frac{D(F_n\circ  v)}{D(F_n)}$ can be  arbitrarily large for sufficiently large $n$. 

(ii)   If $G=K_{1,n}$ $(n\geqslant 6)$ and $v$ is the central vertex of $G$, then $G\circ v=K_n$, and so $D'(G\circ v)=2$. Thus $lim_{n\rightarrow \infty} \frac{D'(G\circ  v)}{D'(G)}=0$.\qed

\medskip

Now we  examine the effects on $D(G)$ and $D'(G)$ when $G$ is modified by an edge contraction. In a graph $G$, contraction of an edge $e$ with endpoints $u,v$ is the replacement of $u$ and $v$ with a single vertex such that edges incident to the new vertex are the edges other than $e$ that were incident with $u$ or $v$. The resulting graph $G\circ e$ has one less edge than $G$ (\cite{Bondy}).
The following theorem gives  a lower bound for the distinguishing number and the distinguishing  index of $G\circ e$:
\begin{theorem}\label{edgecontraction}
Let $G$ be a connected graph of order $n\geqslant 2$ and $e\in E(G)$. We have 
\begin{itemize}
\item[(i)] $D(G)-1 \leqslant D(G\circ e) \leqslant 3D(G)$.
\item[(ii)] $D'(G)-1 \leqslant D'(G\circ e) \leqslant 3D'(G)$.
\end{itemize}
\end{theorem}
\proof
 Let $e$ be the edge  between the two vertices $v$ and $w$. The contraction of $e$, convert $v$ and $w$ to one vertex, and we denote it by $v$, again. 
 
 (i) First we prove $D(G)-1 \leqslant D(G\circ e)$. We define a  distinguishing vertex labeling for $G$ with $D(G\circ e)+1$ labels. For this purpose,  we label $G\circ e$ in a distinguishing way with $D(G\circ e)$ labels. Next we consider  this labeling as a labeling of $G$ by assigning a new label to the vertex $w$. This labeling is a distinguishing labeling for $G$, because if $f$  is an automorphism of $G$ preserving the labeling then $f(w)=w$, and so the restriction of $f$ to $G\circ e$ is an automorphism of $G\circ e$. Since we labeled $G\circ e$ in a distinguishing way at first,  $f$ is the identity automorphism.
 
 Now we prove that $D(G\circ e) \leqslant 3D(G)$. Let $x$ be  a vertex of $G$ and  $N_1(x)$ be the set of vertices of degree one in $G$ which  are adjacent to $x$. Also we use the notation $N(x,y)$, $x,y\in V(G)$, for the set of vertices of $G$ which  are adjacent to both of $x$ and $y$. Using these notations we define a distinguishing vertex labeling of $G\circ e$ with at most $3D(G)$ labels.   We label $G$ with the labels $\{1,\ldots , D(G)\}$ in a distinguishing way. Now we add the number $D(G)$ to the label of each vertex in $N_1(v)$, and next we add the number $2D(G)$ to the label of each vertex in $N(v,w)$. Transfering this labeling to $G\circ e$,  we have a distinguishing labeling for $G\circ e$, because if $f$ is an automorphism of $G\circ e$ preserving the labeling , then $f(N_1(v))=N_1(v)$   and $f(N(v,w))=N(v,w)$. Hence by defining $f(w)=w$, we can extend $f$ to an automorphism of $G$ preserving the labeling. Since we labeled $G$ in a distinguishing way at first, $f$ is the identity automorphism.

  (ii) Here we prove that $D'(G)-1 \leqslant D'(G\circ e)$. We define a  distinguishing edge labeling for $G$ with $D'(G\circ e)+1$ labels. First we label the edge set of $G\circ e$ in a distinguishing way with $D'(G\circ e)$ labels. Assigning a new label to the edge $e$,  we transfer this labeling to $G$. This labeling is a distinguishing labeling for $G$, because if $f$  is an automorphism of $G$ preserving the labeling then $f(\{v,w\})=\{v,w\}$. So the restriction of $f$ to $G\circ e$  is an automorphism of $G\circ e$. Since we labeled $G\circ e$ in a distinguishing way at first,  $f$ is the identity automorphism on $G\circ e$. On the other hand $f$ preserves the adjacency relation on $G$, and so $f(v)=v$ and $f(w)=w$. Therefore  $f$ is the identity automorphism on $G$.
 
  Finally we prove that $D'(G\circ e) \leqslant 3D'(G)$. Let $x\in V(G)$ and $E_1(x)$ be  the set of edges of $G$ that are incident to $x$ and a vertex of degree one of $G$. The set of edges of $G$ that are incident  to $x$ and a vertex of $N(x,y)$, or the edges of $G$ that are incident to $y$ and a vertex of $N(x,y)$ are denoted by $E(x,y)$. Using these notations we define a distinguishing edge labeling of $G\circ e$ with at most $3D'(G)$ labels.   First we label the edge set of $G$ with the labels $\{1,\ldots , D'(G)\}$ in a distinguishing way. Now we add the number $D'(G)$ to the label of each edge in $E_1(v)$, and next we add the number $2D'(G)$ to the label of each edge in $E(v,w)$. It is clear that this labeling is a distinguishing edge labeling of $G$. We want to show that this labeling is a distinguishing edge labeling  for $G\circ e$, too. For this purpose, suppose that  $f$ is an automorphism of $G\circ e$ preserving the labeling, then $f(E_1(v))=E_1(v)$   and $f(E(v,w))=E(v,w)$ (note that the size of $E(v,w)$ decrease to half, after the contraction of $e$). Hence by defining $f(w)=w$, we can extend $f$ to an automorphism of $G$ preserving the labeling. Since we labeled $G$ in a distinguishing way at first, $f$ is the identity automorphism.\qed

\begin{figure}[ht]
\begin{center}
\includegraphics[width=0.57\textwidth]{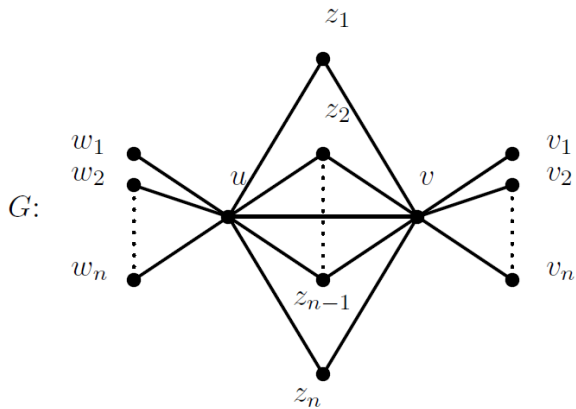}
	\caption{\label{fig8} $D(G\circ e)=3D(G)$.}
\end{center}
\end{figure}

The bounds of Theorem \ref{edgecontraction} for the distinguishing number and index are sharp. Let $K_{1,n}$ be the star graph with $n\geqslant 3$ vertices of degree one. It can be seen that $D(K_{1,n}\circ e)=n-1=D(K_{1,n})-1$ where $e\in E(K_{1,n})$ is an arbitrary edge of $K_{1,n}$. For the second inequality, let $G$ be a graph of order $2n+3$ as shown in Figure \ref{fig8} and $e=uv$. In regard to the degree sequences of $G$ we can obtain the automorphism group, and so we can get $D(G)=n$ and $D(G\circ e)=D(K_{1,3n})=3n$. Theses two  examples are satisfied for showing the sharpness of the distinguishing index in the Theorem \ref{edgecontraction}.

\medskip

\end{document}